\documentclass[reqno]{amsart}
\usepackage{amssymb}
\usepackage[colorlinks=false,hypertex]{hyperref}
\allowdisplaybreaks[4]
\theoremstyle{plain}
\newtheorem{thm}{Theorem}

\theoremstyle{remark}
\newtheorem{rem}{Remark}
\DeclareMathOperator{\td}{d}
\date{Completed on 7 June 2009}
\date{}

\begin{document}

\title[A completely monotonic function involving polygamma functions]
{A completely monotonic function involving the tri- and tetra-gamma functions}

\author[F. Qi]{Feng Qi}
\address[F. Qi]{Department of Mathematics, College of Science, Tianjin Polytechnic University, Tianjin City, 300160, China}
\email{\href{mailto: F. Qi <qifeng618@gmail.com>}{qifeng618@gmail.com}, \href{mailto: F. Qi <qifeng618@hotmail.com>}{qifeng618@hotmail.com}, \href{mailto: F. Qi <qifeng618@qq.com>}{qifeng618@qq.com}}
\urladdr{\url{http://qifeng618.spaces.live.com}}

\author[B.-N. Guo]{Bai-Ni Guo}
\address[B.-N. Guo]{School of Mathematics and Informatics, Henan Polytechnic University, Jiaozuo City, Henan Province, 454010, China}
\email{\href{mailto: B.-N. Guo <bai.ni.guo@gmail.com>}{bai.ni.guo@gmail.com}, \href{mailto: B.-N. Guo <bai.ni.guo@hotmail.com>}{bai.ni.guo@hotmail.com}}

\begin{abstract}
The psi function $\psi(x)$ is defined by $\psi(x)=\frac{\Gamma'(x)}{\Gamma(x)}$ and $\psi^{(i)}(x)$ for $i\in\mathbb{N}$ denote the polygamma functions, where $\Gamma(x)$ is the gamma function. In this paper we prove that a function involving the difference between $[\psi'(x)]^2+\psi''(x)$ and a proper fraction of $x$ is completely monotonic on $(0,\infty)$.
\end{abstract}



\thanks{The authors were supported partially by the Science Foundation of Tianjin Polytechnic University}

\thanks{This paper was typeset using \AmS-\LaTeX}

\maketitle

\section{Introduction}

We recall~\cite[Chapter~XIII]{mpf-1993} and~\cite[Chapter~IV]{widder} that a function
$f$ is said to be completely monotonic on an interval $I$ if $f$ has derivatives of all
orders on $I$ and
\begin{equation}\label{CM-dfn}
(-1)^{n}f^{(n)}(x)\ge0
\end{equation}
for $x\in I$ and $n\ge0$.
\par
The famous Bernstein-Widder Theorem (see \cite[p.~160, Theorem~12a]{widder}) states that a function $f(x)$ on $[0,\infty)$ is completely monotonic if and only if there exists a bounded and non-decreasing function $\alpha(t)$ such that
\begin{equation} \label{berstein-1}
f(x)=\int_0^\infty e^{-xt}\td\alpha(t)
\end{equation}
converges for $x\in[0,\infty)$. This says that a completely monotonic function $f(x)$ on $[0,\infty)$ is a Laplace transform of the measure $\alpha(t)$.
\par
We also recall that the classical Euler gamma function $\Gamma(x)$ is defined by
\begin{equation}\label{gamma-dfn}
\Gamma(x)=\int^\infty_0t^{x-1} e^{-t}\td t,\quad x>0.
\end{equation}
The logarithmic derivative of $\Gamma(x)$, denoted by $\psi(x)=\frac{\Gamma'(x)}{\Gamma(x)}$, is called the psi or di-gamma function, and the derivatives $\psi^{(i)}(x)$ for $i\in\mathbb{N}$ are respectively called the polygamma functions. In particular, the functions $\psi'(x)$ and $\psi''(x)$ are called the tri-gamma and tetra-gamma functions.
\par
In~\cite[p.~208, (4.39)]{forum-alzer}, it was established that the inequality
\begin{equation}\label{psi'psi''}
[\psi'(x)]^2+\psi''(x)>\frac{p(x)}{900x^4(x+1)^{10}}
\end{equation}
holds for $x>0$, where
\begin{equation}
\begin{split}\label{p(x)-dfn}
p(x)&=75x^{10}+900x^9+4840x^8+15370x^7+31865x^6+45050x^5\\
&\quad+44101x^4+29700x^3+13290x^2+3600x+450.
\end{split}
\end{equation}
\par
The aim of this paper is to prove the complete monotonicity of the difference between two functions on both sides of the inequality~\eqref{psi'psi''}.
\par
Our main result may be stated as the following theorem.

\begin{thm}\label{x-4-upper-thm-2}
The function
\begin{equation}\label{x-4-upper-g(x)}
g(x)=[\psi'(x)]^2+\psi''(x)-\frac{p(x)}{900x^4(x+1)^{10}}
\end{equation}
is completely monotonic on $(0,\infty)$, where the function $p(x)$ is defined by~\eqref{p(x)-dfn}.
\end{thm}

\begin{rem}
By the definition of completely monotonic functions, it is easy to see that Theorem~\ref{x-4-upper-thm-2} is stronger than the inequality~\eqref{psi'psi''}, so Theorem~\ref{x-4-upper-thm-2} generalizes the inequality~\eqref{psi'psi''}.
\end{rem}

\section{Proof of Theorem~\ref{x-4-upper-thm-2}}

By the recursion formula \cite[pp.~258 and 260, 6.3.5 and 6.4.6]{abram}
\begin{equation}\label{psisymp4}
\psi^{(n-1)}(x+1)=\psi^{(n-1)}(x)+\frac{(-1)^{n-1}(n-1)!}{x^n}
\end{equation}
for $x>0$ and $n\in\mathbb{N}$, direct calculation produces
\begin{gather*}
g(x)-g(x+1)=\bigl[\psi'(x)-\psi'(x+1)\bigr]\bigl[\psi'(x)+\psi'(x+1)\bigr]\\
+\bigl[\psi''(x)-\psi''(x+1)\bigr] -\biggl[\frac{p(x)}{900x^4(x+1)^{10}} -\frac{p(x+1)}{900(x+1)^4(x+2)^{10}}\biggr]
\end{gather*}
\begin{align*}
&=\frac1{x^2}\biggl[2\psi'(x)-\frac1{x^2}\biggr]-\frac{2}{x^3}
-\biggl[\frac{p(x)}{900x^4(x+1)^{10}} -\frac{p(x+1)}{900(x+1)^4(x+2)^{10}}\biggr]\\
&=\frac2{x^2}\biggl[\psi'(x)-\frac1{2x^2}-\frac1{x} -\frac{p(x)}{1800x^2(x+1)^{10}} +\frac{x^2 p(x+1)}{1800(x+1)^4(x+2)^{10}}\biggr]\\
&=\frac2{x^2}\biggl[\psi'(x)-\frac{1}{2 x}-\frac{3}{4 x^2}-\frac{28}{5 (x+1)}+\frac{51}{10
   (x+2)}+\frac{251}{120 (x+1)^2}\\
&\quad+\frac{331}{120(x+2)^2}-\frac{7}{6 (x+1)^3}+\frac{17}{12 (x+2)^3}+\frac{13}{90
   (x+1)^4}+\frac{49}{72 (x+2)^4}\\
&\quad-\frac{13}{180(x+1)^5}+\frac{47}{180 (x+2)^5}-\frac{1}{120
   (x+1)^6}-\frac{1}{60 (x+2)^6}\\
&\quad+\frac{1}{180 (x+1)^7}-\frac{2}{45(x+2)^7}+\frac{1}{200 (x+1)^8}-\frac{13}{600(x+2)^8}\\
&\quad+\frac{1}{900 (x+1)^9}-\frac{1}{450(x+2)^9}-\frac{1}{1800 (x+1)^{10}}+\frac{1}{450
   (x+2)^{10}}\biggr]\\
&\triangleq \frac2{x^2}H(x).
\end{align*}
Using the formula \cite[p.~255, 6.1.1]{abram}
\begin{equation}\label{gam}
\frac1{x^r}=\frac1{\Gamma(r)}\int_0^\infty t^{r-1}e^{-xt}\td t
\end{equation}
for $r>0$ and $x>0$ and the integral representations \cite[p.~260, 6.4.1]{abram}
\begin{equation} \label{psin}
\psi ^{(n)}(x)=(-1)^{n+1}\int_{0}^{\infty}\frac{t^{n}}{1-e^{-t}}e^{-xt}\td t
\end{equation}
for $n\in\mathbb{N}$ and $x>0$ gives
\begin{align*}
H(x)&=\int_0^\infty\biggl(\frac{t}{1-e^{-t}}-\frac12-\frac34t -\frac{28}5e^{-t}+\frac{51}{10}e^{-2t} +\frac{251}{120}te^{-t}\\
&\quad+\frac{331}{120}te^{-2t}-\frac{7}{12}t^2e^{-t}+\frac{17}{24}t^2e^{-2t} +\frac{13}{540}t^3e^{-t}+\frac{49}{432}t^3e^{-2t}\\
&\quad-\frac{13}{4320}t^4e^{-t} +\frac{47}{4320}t^4e^{-2t}-\frac1{14400}t^5e^{-t} -\frac1{7200}t^5e^{-2t}\\
&\quad +\frac1{129600}t^6e^{-t}-\frac2{32400}t^6e^{-2t}+\frac1{1008000}t^7e^{-t} -\frac{13}{3024000}t^7e^{-2t}\\
&\quad+\frac1{36288000}t^8e^{-t} -\frac1{18144000}t^8e^{-2t} -\frac1{653184000}t^9e^{-t}\\
&\quad+\frac1{163296000}t^9e^{-2t}\biggr)e^{-xt}\td t\\
&=\frac1{653184000}\int_0^\infty\frac{1}{e^t-1}\bigl[163296000 (t-2)e^{3t} -e^{2 t} \bigl(t^9-18 t^8\\ &\quad-648 t^7-5040 t^6+45360 t^5+1965600t^4-15724800 t^3\\
&\quad+381024000 t^2-1856131200 t+3331238400\bigr)+e^t \bigl(5t^9-54 t^8\\
&\quad-3456 t^7-45360 t^6-45360 t^5+9072000 t^4+58363200t^3\\
&\quad+843696000 t^2+435456000 t+6989068800\bigr)-4\bigl(t^9-9 t^8\\
&\quad-702 t^7-10080 t^6-22680 t^5+1776600 t^4+18522000t^3\\
&\quad+115668000 t^2+450424800 t+832809600\bigr)\bigr]e^{-(x+2)t}\td t\\
&\triangleq \frac1{653184000}\int_0^\infty\frac{1}{e^t-1} \theta(t)e^{-(x+2)t}\td t.
\end{align*}
Straightforward computation yields
\begin{align*}
\theta'(t)&=163296000  (3 t-5)e^{3 t}-\bigl(4806345600-2950214400t\\
&\quad+714873600 t^2-23587200 t^3+4158000 t^4+60480 t^5-14616 t^6\\
&\quad-1440 t^7-27 t^8+2 t^9\bigr)e^{2 t} + \bigl(7424524800+2122848000 t\\
&\quad+1018785600 t^2+94651200 t^3+8845200 t^4-317520 t^5-69552 t^6\\
&\quad-3888 t^7-9 t^8+5 t^9\bigr)e^t-36 \bigl(50047200+25704000 t+6174000 t^2\\
&\quad+789600t^3-12600 t^4-6720 t^5-546 t^6-8 t^7+t^8\bigr),\\
\theta''(t)&=489888000 (3 t-4) e^{3 t}-4 \bigl(1665619200-1117670400t\\
&\quad+339746400 t^2-7635600 t^3+2154600 t^4+8316 t^5-9828 t^6\\
&\quad-774 t^7-9 t^8+t^9\bigr) e^{2 t}+ \bigl(9547372800+4160419200 t\\
&\quad+1302739200 t^2+130032000 t^3+7257600 t^4-734832 t^5\\
&\quad-96768 t^6-3960 t^7+36 t^8+5 t^9\bigr)e^t-144 \bigl(6426000+3087000 t\\
&\quad+592200 t^2-12600t^3-8400 t^4-819 t^5-14 t^6+2 t^7\bigr),\\
\theta^{(3)}(t)&=4408992000  (t-1)e^{3 t}-4 \bigl(2213568000-1555848000t\\
&\quad+656586000 t^2-6652800 t^3+4350780 t^4-42336 t^5-25074 t^6\\
&\quad-1620 t^7-9 t^8+2 t^9\bigr) e^{2 t}+ \bigl(13707792000+6765897600 t\\
&\quad+1692835200 t^2+159062400 t^3+3583440 t^4-1315440 t^5\\
&\quad-124488 t^6-3672 t^7+81 t^8+5 t^9\bigr)e^t-1008 \bigl(441000+169200 t\\
&\quad-5400 t^2-4800t^3-585 t^4-12 t^5+2 t^6\bigr),\\
\theta^{(4)}(t)&=4408992000  (3 t-2)e^{3 t}-16 \bigl(717822000-449631000t\\
&\quad+323303400 t^2+1024380 t^3+2122470 t^4-58779 t^5-15372 t^6\\
&\quad-828 t^7+t^9\bigr)e^{2 t} + \bigl(20473689600+10151568000 t+2170022400 t^2\\
&\quad+173396160 t^3-2993760 t^4-2062368 t^5-150192 t^6-3024 t^7\\
&\quad+126 t^8+5 t^9\bigr)e^t-12096 \bigl(14100-900 t-1200 t^2-195 t^3-5
t^4+t^5\bigr),\\
\theta^{(5)}(t)&=13226976000  (3 t-1)e^{3 t}-16 \bigl(986013000-252655200t\\
&\quad+649679940 t^2+10538640 t^3+3951045 t^4-209790 t^5\\
&\quad-36540 t^6-1656 t^7+9 t^8+2 t^9\bigr) e^{2 t}+ \bigl(30625257600\\
&\quad+14491612800 t+2690210880 t^2+161421120 t^3-13305600 t^4\\
&\quad-2963520 t^5-171360 t^6-2016 t^7+171 t^8+5 t^9\bigr)e^t\\
&\quad-60480 \bigl(-180-480 t-117 t^2-4
t^3+t^4\bigr),\\
\theta^{(6)}(t)&=119042784000t e^{3 t}-64 \bigl(429842700+198512370t\\
&\quad+332743950 t^2+9220365 t^3+1713285 t^4-159705 t^5\\
&\quad-21168 t^6-810 t^7+9 t^8+t^9\bigr) e^{2 t}+ \bigl(45116870400\\
&\quad+19872034560 t+3174474240 t^2+108198720 t^3-28123200 t^4\\
&\quad-3991680 t^5-185472 t^6-648 t^7+216 t^8+5 t^9\bigr)e^t\\
&\quad-120960 \bigl(-240-117 t-6 t^2+2 t^3\bigr),\\
\theta^{(7)}(t)&=119042784000  (1+3 t)e^{3 t}-64
\bigl(1058197770+1062512640t\\
&\quad+693148995 t^2+25293870 t^3+2628045 t^4-446418 t^5-48006 t^6\\
&\quad-1548 t^7+27 t^8+2 t^9\bigr)e^{2 t}+ \bigl(64988904960+26220983040 t\\
&\quad+3499070400 t^2-4294080 t^3-48081600 t^4-5104512 t^5\\
&\quad-190008 t^6+1080 t^7+261 t^8+5
t^9\bigr)e^t-362880 \bigl(2 t^2-39-4 t\bigr),\\
\theta^{(8)}(t)&=357128352000  (2+3 t)e^{3 t}-128 \bigl(1589454090+1755661635t\\
&\quad+731089800 t^2+30549960 t^3+1512000 t^4-590436 t^5-53424 t^6\\
&\quad-1440 t^7+36 t^8+2 t^9\bigr) e^{2 t}+ \bigl(91209888000+33219123840 t\\
&\quad+3486188160 t^2-196620480 t^3-73604160 t^4-6244560 t^5\\
&\quad-182448 t^6+3168 t^7+306 t^8+5 t^9\bigr)e^t
-1451520(t-1),\\
\theta^{(9)}(t)&=3214155168000  (1+t)e^{3 t}-128  \bigl(4934569815+4973502870t\\
&\quad+1553829480 t^2+67147920 t^3+71820 t^4-1501416 t^5-116928 t^6\\
&\quad-2592 t^7+90 t^8+4 t^9\bigr)e^{2 t}+ \bigl(124429011840+40191500160 t\\
&\quad+2896326720t^2-491037120 t^3-104826960 t^4-7339248 t^5\\
&\quad-160272 t^6+5616 t^7+351 t^8+5 t^9\bigr)e^t-1451520,\\
\theta^{(10)}(t)&=e^t \bigl[164620512000+45984153600 t+1423215360 t^2-910344960 t^3\\
&\quad-141523200 t^4-8300880 t^5-120960 t^6+8424 t^7+396 t^8+5 t^9\\
&\quad-512\bigl(3710660625+3263666175 t+827275680 t^2+33645780 t^3\\
&\quad-1840860t^4-926100 t^5-63000 t^6-1116 t^7+54 t^8+2 t^9\bigr)e^t\\
&\quad +3214155168000(4+3t)e^{2t}\bigr]\\
&\triangleq e^t\theta_1(t),\\
\theta_1'(t)&=3214155168000 (11+6 t)e^{2 t} -512 \bigl(6974326800+4918217535t\\
&\quad+928213020 t^2+26282340 t^3-6471360 t^4-1304100 t^5-70812 t^6\\
&\quad-684 t^7+72 t^8+2 t^9\bigr) e^t+9 \bigl(5109350400+316270080 t\\
&\quad-303448320
t^2-62899200 t^3-4611600 t^4-80640 t^5\\
&\quad+6552 t^6+352 t^7+5 t^8\bigr),\\
\theta_1''(t)&=8 \bigl[1607077584000  (7+3 t)e^{2 t}-64
\bigl(11892544335+6774643575 t\\
&\quad +1007060040 t^2+396900 t^3-12991860 t^4-1728972 t^5-75600 t^6\\
&\quad-108 t^7 +90 t^8+2 t^9\bigr)e^t +9 \bigl(39533760-75862080
t-23587200 t^2\\
&\quad-2305800 t^3-50400 t^4 +4914 t^5+308 t^6+5 t^7\bigr)\bigr],\\
\theta_1^{(3)}(t)&=8 \bigl[1607077584000 (17+6 t)e^{2 t} -64
\bigl(18667187910+8788763655 t\\
&\quad+1008250740 t^2-51570540 t^3-21636720 t^4-2182572 t^5\\
&\quad-76356 t^6+612 t^7+108 t^8+2 t^9\bigr) e^t-63 \bigl(10837440+6739200t\\
&\quad+988200 t^2+28800 t^3-3510 t^4-264 t^5-5 t^6\bigr)\bigr],\\
\theta_1^{(4)}(t)&=16 \bigl[3214155168000 (10+3 t) e^{2 t}-32
\bigl(27455951565\\
&\quad+10805265135 t+853539120 t^2-138117420 t^3-32549580 t^4\\
&\quad-2640708 t^5-72072 t^6+1476 t^7+126 t^8+2 t^9\bigr)e^t\\
&\quad-945 \bigl(224640+65880t+2880 t^2-468 t^3-44 t^4-t^5\bigr)\bigr],\\
\theta_1^{(5)}(t)&=16 \bigl[3214155168000 (23+6 t)e^{2 t} -32
\bigl(38261216700\\
&\quad+12512343375 t+439186860 t^2-268315740 t^3-45753120 t^4\\
&\quad-3073140 t^5-61740 t^6+2484 t^7+144 t^8+2 t^9\bigr) e^t\\
&\quad-945 \bigl(65880+5760 t-1404
t^2-176 t^3-5 t^4\bigr)\bigr],\\
\theta_1^{(6)}(t)&=64 \bigl[3214155168000 (13+3 t) e^{2 t}-8
\bigl(50773560075\\
&\quad+13390717095 t-365760360 t^2-451328220 t^3-61118820 t^4\\
&\quad-3443580 t^5-44352 t^6+3636 t^7+162 t^8+2 t^9\bigr)e^t\\
&\quad-945 \bigl(1440-702 t-132
t^2-5 t^3\bigr)\bigr],\\
\theta_1^{(7)}(t)&=64 \bigl[3214155168000 (29+6 t) e^{2 t}-8
\bigl(64164277170\\
&\quad+12659196375 t-1719745020 t^2-695803500 t^3-78336720 t^4\\
&\quad-3709692 t^5-18900 t^6+4932 t^7+180 t^8+2 t^9\bigr) e^t\\
&\quad+2835 \bigl(234+88 t+5
t^2\bigr)\bigr],\\
\theta_1^{(8)}(t)&=128 \bigl[6428310336000 (16+3 t)e^{2 t} -4
\bigl(76823473545\\
&\quad+9219706335 t-3807155520 t^2-1009150380 t^3-96885180 t^4\\
&\quad-3823092
t^5+15624 t^6+6372 t^7+198 t^8+2 t^9\bigr)e^t+2835 (44+5 t)\bigr],\\
\theta_1^{(9)}(t)&=128 \bigl[6428310336000  (35+6 t)e^{2 t}-4
\bigl(86043179880\\
&\quad+1605395295 t-6834606660 t^2-1396691100 t^3-116000640 t^4\\
&\quad-3729348
t^5+60228 t^6+7956 t^7+216 t^8+2 t^9\bigr)e^t+14175\bigr],\\
\theta_1^{(10)}(t)&=512 e^t \bigl[6428310336000 (19+3
t)e^t -87648575175+12063818025 t\\
&\quad+11024679960 t^2+1860693660 t^3+134647380 t^4+3367980 t^5\\
&\quad-115920 t^6-9684 t^7-234 t^8-2 t^9\bigr]\\
&\triangleq512 e^t\theta_2(t),\\
\theta_2'(t)&=9 \bigl[1340424225+2449928880 t+620231220 t^2+59843280 t^3\\
&\quad+1871100
t^4-77280 t^5-7532 t^6-208 t^7-2 t^8\\
&\quad+714256704000  (22+3 t)e^t\bigr],\\
\theta_2''(t)&=72 \bigl[306241110+155057805 t+22441230 t^2+935550 t^3\\
&\quad-48300 t^4-5649
t^5-182 t^6-2 t^7+89282088000 (25+3 t)e^t \bigr],\\
\theta_2^{(3)}(t)&=504 \bigl[22151115 + 6411780 t + 400950 t^2 - 27600 t^3 - 4035 t^4\\
&\quad -
   156 t^5 - 2 t^6 + 12754584000 (28 + 3 t)e^t \bigr],\\
\theta_2^{(4)}(t)&=6048\bigl[534315 + 66825 t - 6900 t^2 - 1345 t^3 - 65 t^4 - t^5\\
&\quad +
   1062882000(31 + 3 t) e^t \bigr],\\
\theta_2^{(5)}(t)&=30240 \bigl[13365 - 2760 t - 807 t^2 - 52 t^3 - t^4\\
&\quad +
   212576400(34 + 3 t) e^t \bigr],\\
\theta_2^{(6)}(t)&=60480 \bigl[106288200(37 + 3 t) e^t -1380 - 807 t - 78 t^2 - 2 t^3\bigr],\\
\theta_2^{(7)}(t)&=181440 \bigl[35429400 (40+3 t)e^t -269-52 t-2 t^2\bigr],\\
\theta_2^{(8)}(t)&=725760 \bigl[8857350 (43+3 t) e^t-13-t\bigr],\\
\theta_2^{(9)}(t)&=725760 \bigl[8857350 (46+3 t) e^t-1\bigr].
\end{align*}
It is easy to calculate that
\begin{align*}
  \theta'(0)&=0, & \theta''(0)&=0,\\
  \theta^{(3)}(0)&=0, & \theta^{(4)}(0)&=0,\\
  \theta^{(5)}(0)&=1632960000, & \theta^{(6)}(0)&=17635968000,\\
  \theta^{(7)}(0)&=116321184000, & \theta^{(8)}(0)&=602017920000,\\
  \theta^{(9)}(0)&=2706957792000, & \theta^{(10)}(0)&=11121382944000,\\
  \theta_1'(0)&=31830835680000, & \theta_1''(0)&=83910208435200,\\
  \theta_1^{(3)}(0)&=208999489144320, & \theta_1^{(4)}(0)&=500203983121920,\\
  \theta_1^{(5)}(0)&=1163218362768000, & \theta_1^{(6)}(0)&=2648180949926400, \\ \theta_1^{(7)}(0)&=5932619924353920, & \theta_1^{(8)}(0)&=13125845965639680, \\ \theta_1^{(9)}(0)&=28754776198995840, & \theta_1^{(10)}(0)&=62489726878118400, \\ \theta_2'(0)&=141434891210025, & \theta_2''(0)&=160729807759920,\\
  \theta_2^{(3)}(0)&=180003853569960, &  \theta_2^{(4)}(0)&=199280851953120, \\ \theta_2^{(5)}(0)&=218562955581600, &  \theta_2^{(6)}(0)&=237847398969600, \\ \theta_2^{(7)}(0)&=257132364632640, &   \theta_2^{(8)}(0)&=276417335013120, \\ \theta_2^{(9)}(0)&=295702274730240. &
\end{align*}
Since $\theta_2^{(9)}(t)$ is increasing, so $\theta_2^{(9)}(t)>0$ on $(0,\infty)$,
which means that $\theta_2^{(8)}(t)$ is increasing and positive on $(0,\infty)$. By the
same argument, it is derived that the functions $\theta_2^{(i)}(t)$ for $1\le i\le9$,
$\theta_1^{(i)}(t)$ and $\theta^{(i)}(t)$ for $1\le i\le10$ are increasing and positive
on $(0,\infty)$. Therefore, the function $\theta(t)$ is increasing and positive on
$(0,\infty)$, which implies that the function $H(x)$ is completely monotonic on
$(0,\infty)$. Because the function $\frac2{x^2}$ is completely monotonic on
$(0,\infty)$ and the product of finite completely monotonic functions are also
completely monotonic, we obtain that the function $g(x)-g(x+1)$ is completely monotonic
on $(0,\infty)$, which is equivalent to
$$
(-1)^k[g(x)-g(x+1)]^{(k)}=(-1)^kg^{(k)}(x)-(-1)^kg^{(k)}(x+1)\ge0
$$
for $k\le0$ on $(0,\infty)$. By induction, we have
\begin{multline*}
(-1)^kg^{(k)}(x)\ge(-1)^kg^{(k)}(x+1)\ge(-1)^kg^{(k)}(x+2)\ge\dotsm\\
\ge(-1)^kg^{(k)}(x+m)\ge\dotsm\ge\lim_{m\to\infty}[(-1)^kg^{(k)}(x+m)]=0
\end{multline*}
for $k\ge0$ on $(0,\infty)$. The proof of Theorem~\ref{x-4-upper-thm-2} is complete.

\begin{rem}
Simplification yields that the function
\begin{equation}
H(x)=\psi'(x)-\frac{Q(x)}{1800 x^2 (1+x)^{10} (2+x)^{10}}
\end{equation}
is completely monotonic on $(0,\infty)$, where
\begin{align*}
    Q(x)&=1382400+21657600 x+162792960 x^2+778137600 x^3\\
&\quad +2645782983 x^4+6789381590 x^5+13626443025 x^6\\
&\quad +21889330810 x^7+28579049475 x^8+30634381522 x^9\\
&\quad +27125436630 x^{10}+19896883200 x^{11}+12088287630 x^{12}\\
&\quad +6063596590 x^{13}+2494770300 x^{14}+832958400 x^{15}\\
&\quad +222060150 x^{16}+46134540 x^{17}+7195500 x^{18}\\
&\quad +792300 x^{19}+54900 x^{20}+1800 x^{21}
\end{align*}
for $x\in(0,\infty)$.
\end{rem}


\begin{thebibliography}{99}

\bibitem{abram}
M. Abramowitz and I. A. Stegun (Eds), \textit{Handbook of Mathematical Functions with
Formulas, Graphs, and Mathematical Tables}, National Bureau of Standards, Applied
Mathematics Series \textbf{55}, 9th printing, Washington, 1970.

\bibitem{forum-alzer}
H. Alzer, \textit{Sharp inequalities for the digamma and polygamma functions}, Forum
Math. \textbf{16} (2004), no.~2, 181\nobreakdash--221.

\bibitem{mpf-1993}
D. S. Mitrinovi\'c, J. E. Pe\v{c}ari\'c and A. M. Fink, \textit{Classical and New
Inequalities in Analysis}, Kluwer Academic Publishers, Dordrecht/Boston/London, 1993.

\bibitem{widder}
D. V. Widder, \textit{The Laplace Transform}, Princeton University Press, Princeton, 1941.

\end{thebibliography}
\end{document}